\documentclass[fleqn]{mat01}
\usepackage{times,mathtimy,amssymb}
\begin{document}

\setcounter{page}{99}
\firstpage{99}

\font\zz=msam10 at 10pt
\def\Box{\mbox{\zz{\char'244}}}

\font\www=mtgub at 10.4pt
\def\piy{\mbox{\www{\char'160}}}
\def\Delt{\mbox{\www{\char'104}}}

\font\bi=tibi at 10.4pt

\newtheorem{theore}{Theorem}
\renewcommand\thetheore{\arabic{section}.\arabic{theore}}
\newtheorem{theor}[theore]{\bf Theorem}
\newtheorem{propo}[theore]{\rm PROPOSITION}
\newtheorem{lem}[theore]{Lemma}
\newtheorem{definit}[theore]{\rm DEFINITION}
\newtheorem{coro}[theore]{\rm COROLLARY}
\newtheorem{rem}[theore]{Remark}
\newtheorem{exampl}[theore]{Example}

\renewcommand{\theequation}{\thesection\arabic{equation}}

\title{On a new unified integral}

\markboth{Mridula Garg and Shweta Mittal}{On a new unified integral}

\author{MRIDULA GARG and SHWETA MITTAL}

\address{Department of Mathematics, University of Rajasthan, Jaipur~302~004, India\\
\noindent E-mail: garghome@satyam.net.in; shwetamittal2003@yahoo.co.in}

\volume{114}

\mon{May}

\parts{2}

\Date{MS received 15 March 2003; revised 27 November 2003}

\begin{abstract}
In the present paper we derive a unified new integral whose integrand
contains products of Fox $H$\!-function and a general class of polynomials
having general arguments. A large number of integrals involving various
simpler functions follow as special cases of this integral.
\end{abstract}


\maketitle

\section{Introduction}

The $H$-function introduced by Fox \cite{1}, will be represented and
defined in the following manner:
\begin{align}
H_{p, q}^{m, n} [x] &= H_{p, q}^{m, n} \left\lbrack x \bigg|
\begin{array}{@{}l@{}}
(a_{1}, \alpha_{1}), \ldots, (a_{p}, \alpha_{p})\\[.2pc]
(b_{1}, \beta_{1}), \ldots, (b_{q}, \beta_{q})
\end{array}\right\rbrack\nonumber\\[.2pc]
&= \frac{1}{2 \pi i} \int_{L} \frac{\prod_{j = 1}^{m} \Gamma (b_{j} -
\beta_{j} \xi) \prod_{j = 1}^{n} \Gamma (1 - a_{j} + \alpha_{j}
\xi)}{\prod_{j = m + 1}^{q} \Gamma (1 - b_{j} + \beta_{j} \xi) \prod_{j
= n + 1}^{p} \Gamma (a_{j} - \alpha_{j} \xi)} x^{\xi} \ {\rm d}\xi.
\end{align}
For the nature of contour $L$ in (1.1), the convergence, existence
conditions and other details of the $H$-function, one can refer to
\cite{3}.

The general class of polynomials introduced by Srivastava \cite{4} is
defined in the following manner:
\begin{equation}
S_{V}^{U} [x] = \sum\limits_{K = 0}^{[V/U]} \frac{(-V)_{UK} A (V, K)}{K!}
x^{K}, \quad V = 0, 1, 2, \ldots,
\end{equation}
where $U$ is an arbitrary positive integer and coefficients $A(V,K), (V,
K \geq 0)$ are arbitrary constants, real or complex.

\setcounter{equation}{0}
\section{Main result}

\begin{align}
&\int_{0}^{\infty} x^{\lambda -1 } [x + a + (x^{2} + 2ax)^{1/2}]^{-\nu}
H_{p, q}^{m, n} [ y \{ x + a + (x^{2} + 2ax)^{1/2} \}^{-\mu}]\nonumber\\
&\qquad\qquad \times S_{V}^{U} [z \{ x + a + (x^{2} + 2ax)^{1/2} \}^{-\alpha}]
\hbox{d}x\nonumber
\end{align}
\begin{align}
&= 2a^{-\nu} \left(\frac{1}{2}a\right)^{\lambda} \Gamma (2 \lambda)
\sum\limits_{K = 0}^{[V/U]} (-V)_{UK} A (V, K)
\frac{(z/a^{\alpha})^{K}}{K!} H_{p + 2, q + 2}^{m, n + 2}\nonumber\\[.2pc]
&\quad\, \times \left\lbrack ya^{-\mu} \bigg|
\begin{array}{@{}l@{}}
(-\nu - \alpha K, \mu), (1 + \lambda - \nu - \alpha K, \mu), (a_{1},
\alpha_{1}), \ldots, (a_{p}, \alpha_{p})\\[.2pc]
(b_{1}, \beta_{1}),\ldots,(b_{q}, \beta_{q}), (1 - \nu - \alpha K, \mu),
(-\nu - \alpha K - \lambda, \mu)
\end{array} \!\right\rbrack,
\end{align}
where
\begin{enumerate}
\renewcommand\labelenumi{(\roman{enumi})}
\leftskip .1pc
\item $\mu > 0, \ \hbox{Re} (\lambda, \nu, \alpha) > 0$,

\item $\hbox{Re} (\lambda) - \hbox{Re} (\nu) - \mu
\mathop{\min}\limits_{1 \leq j \leq m} \ \hbox{Re} \
\left(\dfrac{b_{j}}{\beta_{j}}\right) < 0$.
\end{enumerate}

\begin{proof}

To obtain the result (2.1), we first express Fox $H$-function involved
in its left-hand side in terms of contour integral using eq.~(1.1) and
the general class of polynomials $S_{V}^{U} [x]$ in series form given by
eq.~(1.2). Interchanging the orders of integration and summation (which
is permissible under the conditions stated with (2.1)) and evaluating
the $x$-integral with the help of the result given below \cite{2}:
\begin{align*}
&\int_{0}^{\infty} x^{z - 1} [x + a + (x^{2} + 2ax)^{1/2}]^{-\nu} {\rm
d}x\\
&\quad\ = 2 \nu a^{-\nu} \left(\frac{1}{2}a\right)^{z} [\Gamma (1 + \nu
+ z)]^{-1} \Gamma (2z) \Gamma (\nu - z),\quad 0 < \hbox{Re} (z) < \nu,
\end{align*}
\end{proof}
we easily arrive at the desired result (2.1).

\setcounter{equation}{0}
\section{Special case}

If in the integral (2.1) we reduce $S_{V}^{U} [x]$ to unity and Fox
$H$-function to Gauss hypergeometric function \cite{3}, we arrive at the
following result after a little simplification:
\begin{align}
&\int_{0}^{\infty} x^{\lambda - 1} [x + a + (x^{2} +
2ax)^{1/2}]^{-\nu}\nonumber\\
&\qquad\ \times {_{2}F_{1}} (a, b; c; y (x + a + (x^{2} +
2ax)^{1/2})^{-1}){\rm d}x\nonumber\\
&\quad\ = 2^{1 - \lambda} \nu \Gamma (2 \lambda) a^{\lambda - \nu}
\frac{\Gamma (\nu - \lambda)}{\Gamma (\nu + \lambda + 1)}\nonumber\\
&\qquad\ \times {_{4}F_{3}} (a, b, \nu - \lambda, \nu + 1; c, \nu, \nu +
\lambda + 1; y/a),
\end{align}
where
\begin{equation*}
0 < \hbox{Re} (\lambda) < \hbox{Re} (\nu), |y| < |a|.
\end{equation*}
The importance of the result given by (3.1) lies in the fact that it not
only gives the value of the integral but also `augments' the
coefficients in the series in the integrand to  give a $_{4}F_{3}$ series
as the integrated series.

A number of other integrals involving functions that are special cases
of Fox $H$-function \cite{3} and/or the general class of polynomials
\cite{5} can also be obtained from (2.1) but we do not record them here.

\section*{Acknowledgements}

The authors are thankful to the worthy referee for his very valuable
suggestions. The first author is thankful to the University Grants
Commission, New Delhi for providing necessary financial assistance to
carry out the present work. The authors are thankful to K~C~Gupta,
Jaipur for his useful suggestions.

\end{document}